\documentclass[12pt]{article}
\makeatletter
\usepackage[margin=72.27pt]{geometry}
\usepackage{enumitem}
\SetEnumitemKey{special}{topsep=\z@,itemsep=\z@,
  parsep=\z@,partopsep=\z@}
\setlist{special}
\baselineskip=\the\baselineskip plus 0.2pt minus 0.2pt\relax
\parskip\z@
\usepackage[verbose=silent,factor=700,stretch=14,shrink=14,step=1]{microtype}\pretolerance=-1
\usepackage{amsmath}
\usepackage{amssymb}
\usepackage{amsthm}
\DeclareMathOperator{\BR}{BR}
\newtheorem{lemma}{Lemma}
\newtheorem{theorem}[lemma]{Theorem}

\newtheorem{definition}[lemma]{Definition}

\def\heading#1{\bigskip\medskip
  \noindent\textbf{#1}\par
  \nobreak\vskip\bigskipamount}

\usepackage[authordate,doi=false,footmarkoff]{biblatex-chicago}
\begin{document}

{\LARGE\topskip=0.7in
\centerline{Nash Equilibrium and Axiom of Choice Are Equivalent}\smallskip
\centerline{Conrad Kosowsky%
  \def\thefootnote{*}%
  \footnote{University of Michigan, Department of Economics. Email: coko@umich.edu.}%
  \global\advance\c@footnote\m@ne}
}

\vskip 0.3in

\begin{abstract}
In this paper, I prove that existence of pure-strategy Nash equilibrium in games with infinitely many players is equivalent to the axiom of choice.

\medskip


\medskip

\end{abstract}

The study of Nash equilibrium began with games having finitely many players and only slowly progressed to infinite-player games.\footnote{See \textcite{yang-song-2022} for discussion of these developments.} Recent work by \textcite{yang-song-2022} developed a framework for establishing pure-strategy Nash equilibrium existence with infinitely many players as a consequence of the finite-player case. The authors lay out a general approach that they use to prove equilibrium existence in several specific classes of games, and in this paper, I formalize their approach in a single theorem and prove its equivalence to the axiom of choice.\footnote{For readers without a background in set theory, the axiom of choice states that the Cartesian product of nonempty sets is nonempty. Equivalent statements include Zorn's lemma, the well-ordering principle, Tychonoff's theorem, and the fact that every vector space has a basis.}

My contributions are two-fold. First, proving the equivalence clarifies mathmatical issues around equilibrium existence in infinite-player games. Because infinite-player games incorporate an infinite product of strategy spaces, we require the axiom of choice for such games to be non-degenerate, and a relevant question is how conditions for non-degeneracy relate to conditions for equilibrium existence. \textcite{yang-song-2022} demonstrate that the axiom of choice is sufficient to establish equilibrium existence, and my result shows that it is also necessary. In infinite-player games, we may see phenomena that are unintuitive or otherwise different from finite-player games, and better understanding the conditions for equilibrium existence will be helpful for future research on this topic. See \textcite{voorneveld-2010} and \textcite{rachmilevitch-2016,rachmilevitch-2020} for discussion of this point. For existence of mixed-strategy equilibria in games with infinitely many players, see \textcite{salonen-2010}.
Second, the axiom of choice is an object of interest in its own right, and entire books are devoted to listing equivalent and weaker formulations \autocite{howard-rubin-1998,rubin-rubin-1963,rubin-rubin-1985}. My result is a new equivalence between the axiom of choice and a concept in game theory.

\section{Introduction}

We begin with relevant definitions. Given a topological space $X$, we let $C(X)$ denote the hyperspace containing the nonempty closed subsets of $X$. If $X$ is compact Hausdorff, then $C(X)$ is exactly the nonempty compact subsets of $X$. Throughout this paper, we equip $C(X)$ with the upper Vietoris topology, where we take as a basis the subsets $V\subset C(X)$ such that there exists an open set $U\subset X$ with $V=\{K\in C(X)\colon K\subset U\}$. Under our formalism, a continuous function $f\colon X\longrightarrow C(Y)$ is exactly what in traditional economics parlance we would call a closed-valued, upper-hemicontinuous correspondence from $X$ to $Y$. Recall that two points are \textit{topologially distinguishable} if there exists an open set containing one but not the other, and a topological space is \textit{preregular} if any two topologically distinguishable points can be separated by disjoint open neighborhoods.\footnote{For context, a preregular space is Hausdorff if and only if it is $T_0$.}

The main mathematical object of this paper is a \textit{game}, which is a collection of (nonempty) topological spaces $S_i$ indexed by some (possibly infinite) index set $I$ containing at least two elements. We refer to $S_i$ as player $i$'s \textit{strategy space}, and we often focus on the set of outcomes
\[
S=\prod_iS_i,
\]
where we endow $S$ with the product topology. As is standard, we use $-i$ to denote ``players other than $i$,'' and we write
\begin{align*}
S_{-i}&=\prod_{j\not=i}S_j & s_{-i}&\in S_{-i}
\end{align*}
to refer to strategies of players other than $i$. Informally, a point $s_i\in S_i$ represents one possible choice that player $i$ can make, and a point $s_{-i}$ represents one choice made by every player other than $i$. (Although we need choice to ensure that $S$ and $S_{-i}$ are nonempty, we can still talk about their existence even without choice.) For each player $i$, we define the \textit{best-response correspondence} to be a continuous function $\BR_i\colon S_{-i}\longrightarrow C(S_i)$. The focus of this paper is Nash equilibrium, a standard solution concept in game theory. A \textit{pure-strategy Nash equilibrium} is a point $s\in S$ where for any $i$,
\[
s_i\in\BR_i(s_{-i}),
\]
where $s_i$ is the $i$th coordinate of $s$ and $s_{-i}$ is $s$ with the $i$th coordinate removed.

\section{Results}

Our approach involves using Tychonoff's theorem to show that the graphs of all best-response correspondences have nonemtpy intersection. Throughout this paper, we let $\Gamma(f)$ denote the graph of a function $f$, and for a function $f\colon X\longrightarrow C(Y)$ that maps into a hyperspace, we make no distinction between $\Gamma(f)\subset X\times Y$ and $\Gamma(f)\subset X\times C(Y)$ since these two notions are equivalent. We establish a lemma about the graph of a continuous function into a hyperspace and provide two definitions to specify what class of games we are interested in.

\begin{lemma}
If $Y$ is preregular and compact, then a continuous function $f\colon X\longrightarrow C(Y)$ has a closed graph in $X\times Y$.
\end{lemma}

\begin{proof}
Suppose $(x,y)$ satisfies $y\not\in f(x)$. Because $f(x)$ is a closed subset of $Y$, it is compact, and every point in $f(x)$ is topologically distinguishable from $y$. By compactness of $f(x)$ and preregularity, there exist open sets $U$ and $V$ such that $y\in U$, $f(x)\subset V$, and $U\cap V=\varnothing$. By continuity of $f$, the set $W=\{x\colon f(x)\subset V$\} is open in $X$, so $W\times U\subset X\times Y$ is open, contains $(x,y)$, and does not intersect $\Gamma(f)$. Because $(x,y)$ was arbitrary, it follows that $\Gamma(f)$ is closed.
\end{proof}

\begin{definition} 
Given a set $J\subset I$ containing at least two elements, define
\begin{align*}
S_J&=\prod_{i\in J}S_i & S_{{-}J}&=\prod_{i\in J^{\text c}}S_i
\end{align*}
For a point\vadjust{\penalty-100} $p\in S_{-J}$, define the \textbf{reduction around p} to be the game with (1) index set $J$; (2) strategy spaces $S_i$, $i\in J$; and (3) best-response correspondences given by $\BR_i(p,.)$, $i\in J$. In the event that $J=I$, we may define the (unique) corresponding reduction to be the original game. A reduction is \textbf{finite} if $J$ is finite.
\end{definition}

\begin{definition} A game is \textbf{well-specified} if (1) each $S_i$ is nonempty, preregular, and compact; (2) each best-response correspondence is continuous; and (3) every finite reduction has a pure-strategy Nash equilibrium.
\end{definition}

Definition~3 highlights the importance of equilibrium in finite-player games for the infinite-player case. We come to the main results of the paper. Theorem~4 establishes the existence of a pure-strategy Nash equilibrium, and Theorem~5 proves the axiom of choice assuming that theorem~4 holds. Theorem~4 is very similar to the ideas in \textcite{yang-song-2022}.

\begin{theorem}
Every well-specified game has a pure-strategy Nash equilibrium.
\end{theorem}

\begin{proof}
If $I$ is finite, the result follows immediately from the definition with $J=I$. Consider the case with infinite $I$. From Tychnoff's theorem, we know that $S$ is compact, and Lemma~1 means that each set $\Gamma(\BR_i)^{\text c}$ is open in $S$. Consider any finite set $J\subset I$. The axiom of choice implies that $S_{-J}$ is nonempty, so there exists a finite reduction of the game involving only those $S_i$ with $i\in J$. This finite reduction has a pure-strategy Nash equilibrium, and it follows that
\[
\bigcap_{i\in J}\Gamma(\BR_i)\not=\varnothing.
\]
Because no finite subcollection of $\{\Gamma(\BR_i)^{\text c}\}$ covers $S$, it must be the case that $\{\Gamma(\BR_i)^{\text c}\}$ does not cover $S$. But any point not covered by this collection is a pure-strategy Nash equilibrium.
\end{proof}

\begin{theorem}[Axiom of Choice]
The Cartesian product of nonempty sets is nonempty.
\end{theorem}

\begin{proof}
Let $\{X_i\}$ be a collection of nonempty sets indexed by $I$. If $I$ is finite, the result is true in ZF without choice, so consider the case with $I$ infinite. Endow each set with the indiscrete topology. Then each $X_i$ is compact and preregular, and $C(X_i)$ is a singleton. Thus there exists a unique map $\BR_i\colon X_{-i}\longrightarrow C(X_i)$, and this map must be continuous. If the product of all $X_i$ is empty, then the game has no finite reductions, so trivially every finite reduction has a pure-strategy Nash equilibrium. But then the game has a pure-strategy Nash equilibrium, so the product would be nonempty, a contradiction.
\end{proof}

\section{Conclusion}

We have established that pure-strategy equilibrium existence in infinite games is equivalent to the axiom of choice. Our process for proving one direction relied on Tychonoff's theorem to extend equilibrium from finite reductions to the entire game, and for the other direction, we used a game with indiscrete topological spaces and maximal best responses. The reliance on finite-player games is a natural requirement in that many results from analysis and topology establish pure-strategy Nash equilibrium in finite-player games, and set-theoretic concerns arise only in games with infinitely many players. This result clarifies the mathematical intuition surrounding pure-strategy equilibrium existence in infinite-player games and provides a new formulation of the axiom of choice.

\heading{References}

{
\def\mkbibnamefamily#1{\textsc{#1}}
\def\mkbibnamegiven#1{\textsc{#1}}
\def\mkbibnameprefix#1{\textsc{#1}}
\def\mkbibnamesuffix#1{\textsc{#1}}

\printbibliography[heading=none]

\vskip-\lastskip

}

\end{document}